\documentclass{rmmcart}
\newtheorem{prop}{Proposition}

     \title[On the Area Bounded by $\prod\limits_{k = 1}^n\left|x\sin\left(\frac{k\pi}{n}\right) -y \cos\left(\frac{k\pi}{n}\right)\right| = 1$]{On the Area Bounded by the Curve $\prod\limits_{k = 1}^n\left|x\sin\left(\frac{k\pi}{n}\right) -y\cos\left(\frac{k\pi}{n}\right)\right| = 1$}

     \author{Anton Mosunov}
     \address{University of Waterloo, 200 University Ave W, Waterloo ON, Canada N2L 3G1}
     \email{amosunov@uwaterloo.ca}

     \date{February, 1, 2020}


\begin{document}
\begin{abstract}
For a positive integer $n$, let
$$
F_n^*(X, Y) = \prod\limits_{k = 1}^n\left(X\sin\left(\frac{k\pi}{n}\right) -Y\cos\left(\frac{k\pi}{n}\right)\right).
$$
In 2000 Bean and Laugesen proved that for every $n \geq 3$ the area bounded by the curve $|F_n^*(x, y)| = 1$ is equal to $4^{1 - 1/n}B\left(\frac{1}{2} - \frac{1}{n}, \frac{1}{2}\right)$, where $B(x, y)$ is the beta function. We provide an elementary proof of this fact based on the polar formula for the area calculation. We also prove that
$$
F_n^*(X, Y) = 2^{1 - n}\sum\limits_{\substack{1 \leq k \leq n\\\text{$k$ is odd}}}(-1)^{\frac{k - 1}{2}}\binom{n}{k}X^{n - k}Y^k
$$
and demonstrate that $\ell_n = 2^{n - 1 - \nu_2(n)}$ is the smallest positive integer such that the binary form $\ell_n F_n^*(X, Y)$ has integer coefficients. Here $\nu_2(n)$ denotes the \mbox{$2$-adic} order of $n$.
\end{abstract}
     \maketitle
     
\section{Introduction}

Let $F(X, Y) = a_0X^n + a_1X^{n - 1}Y + \cdots + a_nY^n$ be a binary form of degree $n \geq 3$ with complex coefficients and nonzero discriminant $D_F$. In this article we evaluate the area $A_F$ bounded by the curve
$$
|F(x, y)| = 1, \quad (x, y) \in \mathbb R,
$$
associated with a particular family of binary forms $F$. The quantity $A_F$ arises in the study of \emph{Thue inequalities} $|F(x, y)| \leq h$ where the coefficients of $F$ are integers and $h$ is a fixed positive integer. In 1933 Mahler \cite{mahler33} proved that the number of integer solutions $Z_F(h)$ to the Diophantine inequality above satisfies
$$
\left|Z_F(h) - A_Fh^{2/n}\right| \ll_F h^{1/(n - 1)},
$$
provided that $F$ is irreducible. More recently, Stewart and Xiao \cite{stewart-xiao} proved that the number of integers of absolute value at most $h$ which are represented by $F$ is asymptotic to $C_Fh^{2/n}$, where a positive number $C_F$ depends solely on $F$ and is a rational multiple of $A_F$.

In 1994 Bean \cite{bean94} proved that the quantity $|D_F|^{1/n(n - 1)}A_F$ is bounded by the absolute constant $3B\left(\frac{1}{3}, \frac{1}{3}\right) \approx 15.90$, where
\begin{equation} \label{eq:beta}
B(x, y) = 2\int\limits_0^{\pi/2}(\sin \theta)^{2x - 1}(\cos \theta)^{2y - 1}d\theta
\end{equation}
is the beta function. In his investigation of the quantity
$$
M_n = \max\left\{|D(F)|^{1/n(n - 1)}A_F\right\},
$$
where the maximum is taken over all binary forms $F$ of degree $n$ with nonzero discriminant, Bean \cite{bean95} conjectured that the value of $M_n$ is attainable by the binary form
$$
F_n^*(X, Y) = \prod\limits_{k = 1}^n\left(X\sin\left(\frac{k\pi}{n}\right) -Y \cos\left(\frac{k\pi}{n}\right)\right).
$$
The curve $|F_n^*(x, y)| = 1$ has the special property that it is invariant under rotation by any integer multiple of $\pi/n$. In 2000 Bean and Laugesen \cite{bean-laugesen} proved that for every $n \geq 3$,
\begin{equation} \label{eq:AFnstar}
A_{F_n^*} = 4^{1 - 1/n}B\left(\frac{1}{2} - \frac{1}{n}, \frac{1}{2}\right), \quad D_{F_n^*}^{1/n(n - 1)} = \frac{1}{2}n^{1/(n - 1)}.
\end{equation}
In \mbox{Section \ref{sec:area}} we prove the first of these equalities via elementary methods. Apart from basic trigonometric identities, our proof involves the identity $\sin(n\theta) = 2^{n - 1}\prod_{k = 1}^n\sin\left(\frac{k\pi}{n} - \theta\right)$ (established in \mbox{Section \ref{sec:identity}}), the formula for the area bounded by a curve in polar form; and the trigonometric form of the beta function (\ref{eq:beta}).

In \mbox{Section \ref{sec:coeff}} we explain what the coefficients of $F_n^*(X, Y)$ look like. More precisely, define the binary form
$$
S_n(X, Y) = \ell_nF_n^*(X, Y),
$$
where $\ell_n$ is the smallest positive integer such that $S_n$ has integer coefficients. We prove the following result.

\begin{prop} \label{prop:S}
For any positive integer $n$,
\begin{equation} \label{eq:Fn*}
F_n^*(X, Y) = 2^{1 - n}\sum\limits_{\substack{1 \leq k \leq n\\\text{$k$ is odd}}}(-1)^{\frac{k - 1}{2}}\binom{n}{k}X^{n - k}Y^k
\end{equation}
and
\begin{align} \label{eq:S}
S_n(X, Y)
& = 2^{n - 1 - \nu_2(n)}F_n^*(X, Y)\\\notag
& = 2^{-\nu_2(n)}\sum\limits_{\substack{1 \leq k \leq n\\\text{$k$ is odd}}}(-1)^{\frac{k - 1}{2}}\binom{n}{k}X^{n - k}Y^k,
\end{align}
where $\nu_2(n)$ is the $2$-adic order of $n$.
\end{prop}

Finally, in \mbox{Section \ref{sec:area2}} we derive the formula for $A_{S_n}$:
$$
A_{S_n} = 4^{\nu_2(n)/n}B\left(\frac{1}{2} - \frac{1}{n}, \frac{1}{2}\right).
$$

\section{Proof of the Identity $\sin(n\theta) = 2^{n - 1}\prod_{k = 1}^n\sin\left(\frac{k\pi}{n} - \theta\right)$.} \label{sec:identity}

Let $U_{n - 1}(X)$ be the $(n - 1)$-st Chebyshev polynomial of the second kind:
$$
U_{n - 1}(X) = 2^{n - 1}\prod\limits_{k = 1}^{n - 1}\left(X - \cos\left(\frac{k\pi}{n}\right)\right).
$$
Let $\theta \in \mathbb R$. In view of the identity $\sin(n\theta) = U_{n - 1}(\cos \theta)\sin \theta$, we obtain
\small
\begin{align*}
\sin(n\theta) & = \left[2^{n - 1}\prod\limits_{k = 1}^{n - 1}\left(\cos\theta - \cos\left(\frac{k\pi}{n}\right)\right)\right]\sin \theta\\
& = \left[2^{n - 1}\prod\limits_{k = 1}^{n - 1}\left(2\sin\left(\frac{k\pi}{2n} - \frac{\theta}{2}\right)\sin\left(\frac{k\pi}{2n} + \frac{\theta}{2}\right)\right)\right]\sin \theta\\
& = \left[2^{n - 1}\prod\limits_{k = 1}^{n - 1}\left(2\sin\left(\frac{k\pi}{2n} - \frac{\theta}{2}\right)\right)\prod\limits_{k = 1}^{n - 1}\cos\left(\frac{(n - k)\pi}{2n} - \frac{\theta}{2}\right)\right]\sin \theta\\
& = \left[2^{n - 1}\prod\limits_{k = 1}^{n - 1}\left(2\sin\left(\frac{k\pi}{2n} - \frac{\theta}{2}\right)\cos\left(\frac{k\pi}{2n} - \frac{\theta}{2}\right)\right)\right]\sin(\pi - \theta)\\
& = 2^{n - 1}\prod\limits_{k = 1}^n\sin\left(\frac{k\pi}{n} - \theta\right).
\end{align*}
\normalsize

\section{Proof of the Identity $A_{F_n^*} = 4^{1 - 1/n}B\left(\frac{1}{2} - \frac{1}{n} , \frac{1}{2}\right)$.} \label{sec:area}

Consider the curve $|F_n^*(x, y)| = 1$. The change of variables
$$
x = r(\theta)\cos \theta, \quad y = r(\theta)\sin \theta
$$
enables us to transform this equation into polar coordinates:
\begin{align*}
r(\theta)
& = |F_n^*(\cos \theta, \sin \theta)|^{-1/n}\\
& = \left|\prod\limits_{k = 1}^n\left(\cos\theta\sin\left(\frac{k\pi}{n}\right) - \sin\theta\cos\left(\frac{k\pi}{n}\right)\right)\right|^{-1/n}\\
& = \left|\prod\limits_{k = 1}^n\sin\left(\frac{k\pi}{n} - \theta\right)\right|^{-1/n}\\
& = \left|2^{1 - n}\sin(n\theta)\right|^{-1/n}.
\end{align*}
We are now able to apply a well-known formula for the area bounded by a curve in polar form:
$$
A_{F_n^*} = \frac{1}{2}\int\limits_0^{2\pi}r(\theta)^2d\theta = \frac{1}{2}\int\limits_0^{2\pi}\left|2^{1 - n}\sin(n\theta)\right|^{-2/n}d\theta = \frac{2^{1 - 2/n}}{n}\int\limits_0^{2\pi n}|\sin \theta|^{-2/n}d\theta.
$$
Notice that the function $|\sin \theta|^{-2/n}$ has period $\pi$, and furthermore $\int_0^{\pi/2}\left(\sin \theta\right)^{-2/n}d\theta = \int_{\pi/2}^{\pi}\left(\sin \theta\right)^{-2/n}d\theta$. Hence
\begin{align*}
A_{F_n^*}
& = \frac{2^{1 - 2/n}}{n}\cdot 4n\int\limits_0^{\pi/2}\left(\sin \theta\right)^{-2/n}d\theta\\
& = 4^{1 - 1/n}\cdot 2\int\limits_0^{\pi/2}(\sin \theta)^{-2/n}d\theta\\
& = 4^{1 - 1/n}B\left(\frac{1}{2} - \frac{1}{n}, \frac{1}{2}\right),
\end{align*}
where the last equality follows from (\ref{eq:beta}).

\section{Coefficients of $F_n^*(X, Y)$} \label{sec:coeff}

In this section we prove \mbox{Proposition \ref{prop:S}}. If $n = 1$ the result holds, so we assume that $n \geq 2$. Let $f_n^*(X)$ and $h_n^*(X)$ denote the polynomials
$$
f_n^*(X) = \prod\limits_{k = 1}^{n - 1}\left(\sin\left(\frac{k\pi}{n}\right)X - \cos\left(\frac{k\pi}{n}\right)\right)
$$
and
$$
h_n^*(X) = 2^{1 - n}\operatorname{Im}\left((X + i)^n\right),
$$
where $\operatorname{Im}(z)$ denotes the imaginary part of a complex number $z$. It is a consequence of the Binomial Theorem that
\begin{align*}
h_n^*(X)
& = 2^{1 - n}\operatorname{Im}\left(\sum\limits_{0 \leq k \leq n}\binom{n}{k}X^{n - k}i^k\right)\\
& = 2^{1 - n}\sum\limits_{\substack{1 \leq k \leq n\\\text{$k$ is odd}}}(-1)^{\frac{k - 1}{2}}\binom{n}{k}X^{n - k}.
\end{align*}
We claim that $f_n^* = h_n^*$. To prove this, we will demonstrate that these two polynomials have the same roots, namely $\cot\left(\frac{k\pi}{n}\right)$ for $k = 1, 2, \ldots, n - 1$, as well as equal leading coefficients.

First, notice that both $f_n^*$ and $h_n^*$ have degrees $n - 1$. Second, if $k$ is an integer between $1$ and $n - 1$, then $\sin\left(\frac{k\pi}{n}\right) \neq 0$, and so it immediately follows from the definition of $f_n^*$ that $f_n^*\left(\cot\left(\frac{k\pi}{n}\right)\right) = 0$. In turn, \mbox{De Moivre's} formula tells us that
\begin{align*}
h_n^*\left(\cot\left(\frac{k\pi}{n}\right)\right)
& = 2^{1 - n}\operatorname{Im}\left(\left(\cot\left(\frac{k\pi}{n}\right) + i\right)^n\right)\\
& = 2^{1 - n}\sin\left(\frac{k\pi}{n}\right)^{-n}\operatorname{Im}\left(\left(\cos\left(\frac{k\pi}{n}\right) + i\sin\left(\frac{k\pi}{n}\right)\right)^n\right)\\
& = 2^{1 - n}\sin\left(\frac{k\pi}{n}\right)^{-n}\operatorname{Im}(\cos(k\pi) + i\underbrace{\sin(k\pi)}_{0})\\
& = 0.
\end{align*}
Therefore, polynomials $f_n^*$ and $h_n^*$ have the same roots. Third, notice that the leading coefficient of $f_n^*(X)$ is $\prod_{k = 1}^{n - 1}\sin\left(\frac{k\pi}{n}\right)$, while the leading coefficient of $h_n^*(X)$ is $2^{1 - n}n$. We use the result established in \mbox{Section \ref{sec:identity}} to prove that these coefficients are equal to each other:
$$
\prod_{k = 1}^{n - 1}\sin\left(\frac{k\pi}{n}\right) = \lim\limits_{\theta \rightarrow 0}\prod_{k = 1}^{n - 1}\sin\left(\frac{k\pi}{n} - \theta\right) = \lim\limits_{\theta \rightarrow 0}2^{1 - n}\frac{\sin(n\theta)}{\sin(\theta)} = 2^{1 - n}n.
$$
Therefore, $f_n^* = h_n^*$.

Now, let $f_n^*(X, Y)$ and $h_n^*(X, Y)$ denote the homogenizations of $f_n^*(X)$ and $h_n^*(X)$, respectively, so that $f_n^*(X, Y) = h_n^*(X, Y)$. Then $F_n^*(X, Y) = Yf_n^*(X, Y) = Yh_n^*(X, Y)$, which is the same as (\ref{eq:Fn*}).

It remains to explain why the binary form $S_n = 2^{n - 1 - \nu_2(n)}F_n^*$ has integer coefficients and why the greatest common divisor of its coefficients is equal \mbox{to $1$}. Let
$$
g(n) = \gcd\limits_{\substack{1 \leq k \leq n\\ \text{$k$ is odd}}}\left\{\binom{n}{k}\right\}.
$$
In view of (\ref{eq:S}), we need show that $g(n) = 2^{\nu_2(n)}$.

Assume for a contradiction that there exists an odd prime $p$ such that \mbox{$p \mid g(n)$}. Let $k = p^r$, where $r$ is the largest positive integer such that \mbox{$p^r \mid n$}. Since $k$ is an odd integer between $1$ and $n$, we conclude that $p \mid \binom{n}{k}$. Letting $\nu_p(m)$ denote the $p$-adic order of a positive integer $m$, we conclude that \mbox{$\nu_p\left(\binom{n}{k}\right) > 0$}. One can then use Legendre's formula $\nu_p(m!) = \sum_{j = 1}^\infty \left\lfloor\frac{m}{p^j}\right\rfloor$ to show that $\nu_p\left(\binom{n}{k}\right) = \nu_p(n!) - \nu_p(k!) - \nu_p\left((n - k)!\right) = 0$, which contradicts our observation that the $p$-adic order of $\binom{n}{k}$ is positive. Therefore, $g(n)$ has no odd prime factors, or in other words $g(n) = 2^s$ for some non-negative \mbox{integer $s$}.


It remains to show that $s = \nu_2(n)$. Since $g(n) \mid \binom{n}{1} = n$ and $g(n) = 2^s$, we see that $s \leq \nu_2(n)$. To establish the reverse inequality $s \geq \nu_2(n)$, recall the well-known result
$$
\frac{n}{\gcd(n, k)} \mid \binom{n}{k},
$$
which was attributed by Dickson \cite{dickson} to Hermite. Now, if $k$ is an odd integer between $1$ and $n$, we see that $\gcd(n, k)$ is odd, which means that \mbox{$2^{\nu_2(n)} \mid \frac{n}{\gcd(n, k)}$}. But then $2^{\nu_2(n)} \mid \binom{n}{k}$. Since $k$ was chosen arbitrarily, we conclude that \mbox{$2^{\nu_2(n)} \mid g(n)$}. Since $g(n)= 2^s$, we see that $s \geq \nu_2(n)$, and so we come to the desired conclusion that $s = \nu_2(n)$.

\section{Proof of the Identity $A_{S_n} = 4^{\nu_2(n)/n}B\left(\frac{1}{2} - \frac{1}{n}, \frac{1}{2}\right)$} \label{sec:area2}

The formula $A_{F} = \int_{-\infty}^{+\infty}|F(x, 1)|^{-2/n}dx$ implies that, for any nonzero complex number $c$, $A_{cF} = |c|^{-2/n}A_F$. Combining this result with (\ref{eq:AFnstar}) and (\ref{eq:S}), we find that
\begin{align*}
A_{S_n}
& = A_{2^{n - 1 - \nu_2(n)}F_n^*} = 2^{-\frac{2(n - 1 - \nu_2(n))}{n}}A_{F_n^*} = 4^{\nu_2(n)/n}B\left(\frac{1}{2} - \frac{1}{n}, \frac{1}{2}\right).
\end{align*}

\end{document}